\documentclass[11pt]{article}
	\usepackage{setspace}
\usepackage{tabularx} 
\usepackage{amssymb,amsmath,amsthm, bbm, mathtools} 
\usepackage{float}
\usepackage{enumitem}
\usepackage{authblk}
\usepackage{caption, graphbox, subcaption}
\usepackage{parskip}

\usepackage{tcolorbox}
\usepackage{graphicx, float} 
\usepackage[margin=1in,letterpaper]{geometry} 
\usepackage{cite} 
\usepackage[final]{hyperref} 

\usepackage[colorinlistoftodos,prependcaption,textsize=scriptsize]{todonotes}

\newcommand{\dx}{\mathrm{d}}
    \newtheorem{theorem}{Theorem}
    \newtheorem{lemma}[theorem]{Lemma}
    \newtheorem{proposition}[theorem]{Proposition}
    \newtheorem{corollary}{Corollary}[theorem]
    
    \newtheorem{definition}{Definition}
    \newtheorem{remark}{Remark}

\title{Batch sojourn time in polling systems on a circle}
\date{\today}

\usepackage[square, numbers]{natbib}
\usepackage{tikz}
\usetikzlibrary{decorations.markings}
\usetikzlibrary{shapes.geometric}
\usetikzlibrary{patterns}

\usepackage[section]{placeins}
\begin{document}

\author[1,*]{Tim Engels}
\author[2]{Ivo Adan}
\author[1]{Onno Boxma}
\author[1]{Jacques Resing}

\affil[1]{Department of Mathematics and Computer Science, Eindhoven University of Technology}
\affil[2]{Department of Industrial Engineering and Innovation Sciences, Eindhoven University of Technology}
\affil[*]{Corresponding author: t.p.g.engels@tue.nl}

\maketitle
\begin{abstract}
    In this paper, we analyze a polling system on a circle. Random batches of customers arrive at a circle, where each customer, independently, obtains a uniform location. A single server cyclically travels over the circle to serve all customers. Using mean value analysis, we derive the expected number of waiting customers within a given distance of the server and a closed form expression for the mean batch sojourn time.
\end{abstract}

\section{Introduction}
\label{sec:Intro}
Polling systems are systems of multiple queues, where a single server visits the queues in a cyclic order to complete the service requirements of customers. A vast amount of literature is devoted to the analysis of these systems; see for instance the literature reviews in \citet{Takagi2001, Borst2018}. Polling models naturally arise in a wide range of application areas, for example in communication, traffic and transportation systems \citep{Grillo1990, Levy1990, Boon2011}. Throughout the years, numerous extensions of the model have been analyzed, in particular we focus on two extensions: extension to a continuum of queues, i.e. polling on a circle,
and the case of batch arrivals.\\
Polling on a circle with batch arrivals naturally arises in warehouse logistics. Customer orders, often consisting of several different items, arrive at the warehouse, requesting to be picked in a timely manner. One way of picking these orders, is by using a so-called milk-run system \citep{Gaast2019}. In this system, a picker walks cyclically through the warehouse, picking any requested item that she encounters. This can be modeled as a polling system, where each requested item represents a customer, an order therefore can be seen as a customer batch. The server represents the picker walking through the warehouse. Warehouses often store a very large number of different items and therefore it is natural to approximate the warehouse by allowing for \emph{continuous} pick locations, i.e. extending the model to polling on a circle. The batch sojourn time then represents the time it takes the picker to completely pick an order.\\
Polling on a circle, see e.g. \citet{Coffman1986, Kroese1992, Kroese1993}, assumes that customers arrive according to a Poisson process and obtain locations according to a continuous distribution on the circle. The former two papers assume that these locations are uniform, while the latter considers an extension where the customers arrive according to a certain distribution with respect to the server. \citet{Eliazar2003}, considers an extended model: ``the snowblower problem'', which assumes a more general arrival process of jobs, excluding bulk arrivals, on the circle and allows for a stochastic travel time over intervals. Using an analytical approach, the author derives some results on the amount of work in intervals, in particular the expected amount of work in an interval is found. \citet{Eliazar2005} later proved that these snowblower models are limiting cases of discrete polling systems, where the number of queues tends to infinity. Hence, these models can prove to be useful for the analysis of large and complex polling models. For example, the mean waiting time has an explicit formula for many continuous polling models, whereas under discrete polling models no method has been found that is faster than solving a system of $N^2$ linear equations, where $N$ is the number of queues \citep{Winands2006}. On the other hand obtaining more intricate performance measures than the mean in polling on a circle is more challenging.\\
The extension to batch arrivals assumes that customers now arrive as a batch, where each such group can be considered as a single job. The locations of the customers in these batches are assumed to follow a certain distribution. In the current paper, we assume that these locations are uniformly distributed and independent from one another. The mean waiting time of customers in bulk arrival polling systems is previously analyzed in 
\citet{Levy1991}, providing a mean value analysis. \citet{Boxma1989} derives a conservation law for the weighted mean waiting time of customers in the system, making use of a stochastic decomposition. Other research on this topic includes that of \citet{Mei2002, AlHanbali2012}. The analysis of the batch sojourn time in polling systems is, as far as we know, restricted to discrete systems, see \citet{Gaast2017}. In said paper the authors provide a mean value analysis for the mean sojourn time of a batch, again requiring one to solve a system of $N^2$ linear equations. \\
In the current paper we combine both extensions and present a mean value analysis for the sojourn time of a batch of customers under a continuous variant. For this, we analyze the number of waiting customers in the system within a certain distance of the server. Using Little's law and an explicit expression for the expected waiting time, we derive an integral equation for the number of customers at distance $x$ of the server. Our main contributions are a novel technique used for the analysis of the mean batch sojourn time and an explicit expression for the mean batch sojourn time that lends itself for optimization purposes. The technique in this paper allows for the analysis of related systems, like systems with generally distributed arrival locations of customers. This is a topic for future research.\\
The paper is built up as follows: Section \ref{sec:ModelDescr} is devoted to the description of the model and some preliminaries. We determine the steady-state mean number of waiting customers in Section \ref{sec:MeanCust}. Section \ref{sec:BatchSojourn} is devoted to the determination of the mean of the total sojourn time of a batch of customers. Some numerical results are presented in Section \ref{sec:NumericalResults}.

\section{Model description and preliminaries}
\label{sec:ModelDescr}
Consider a circle with circumference $1$. Assume that batches of customers arrive on this circle according to a Poisson process with intensity $\lambda$; each batch is considered as a single job. Each arriving customer in the  batch is assigned a location on this circle according to the uniform distribution, independent from all other customers in the same batch. The size of a customer batch, $K$, is assumed to be strictly positive and to follow a known distribution with probabilities $p_k$ and probability generating function $\tilde{K}(\cdot)$ (without loss of generality we assume that $p_0 = 0$). \\
The service is provided by a single traveling server. We assume that the server travels at a fixed speed $\alpha^{-1}$ in a given direction and serves any customer she encounters, taking a randomly distributed time $B$. After finishing the service, the server is assumed to continue traveling in the same direction.\\
Extending the methodology of \citet[Theorem 3.1]{Kroese1993} to batch arrivals shows that the system is stable when $\rho := \lambda\mathbb{E}[K]\mathbb{E}[B]<1$, as is intuitively logical. Throughout this paper we consider the limiting behavior of the system and in particular focus on (i) the steady-state number of waiting customers in the system, and (ii) the long-run sojourn time of a job, $S^B$, that is the time from the arrival of the batch of customers to the service completion of the last customer in that batch.

\begin{figure}[!htbp]
\centering
\begin{tikzpicture}
\draw[postaction = {decorate, decoration = {markings, mark = at position 0.6 with {\node[draw, red, circle, fill, inner sep = 1mm] (server){};}}}]
[postaction = {decorate, decoration = {markings, mark = at position 0.8 with {\node[regular polygon, regular polygon sides = 5,draw, black,  fill, inner sep = 0.7mm]{};}}}]
[postaction = {decorate, decoration = {markings, mark = at position 0.35 with {\node[regular polygon, regular polygon sides = 5,draw, black,  fill, inner sep = 0.7mm]{};}}}]
[postaction = {decorate, decoration = {markings, mark = at position 0.55 with {\node[regular polygon, regular polygon sides = 5,draw, black,  fill, inner sep = 0.7mm]{};}}}]
[postaction = {decorate, decoration = {markings,  mark = at position 0.45 with {\node[inner sep = 0mm] (xlabel) {};}}}]
[postaction = {decorate, decoration = {markings,  mark = at position 0.5 with {\arrowreversed[line width = 0.7mm]{stealth}}}}]
[postaction = {decorate, decoration = {markings,  mark = at position 1 with {\arrowreversed[line width = 0.7mm]{stealth}}}}]
(0,0) circle (2);
\node [left] at (server.west) {Server};
\node [above left] at (xlabel) {$x$};
\path[draw] (server)--(0,0);
\filldraw[fill opacity=0.5,fill=blue!30] (0,0) -- (server) arc (216:108:2cm) -- (0,0);
\end{tikzpicture}
\caption{Illustration of the polling model and corresponding range within distance $x$ of the server. In this example $L(x) = 2$.}
\label{fig:model}
\end{figure}
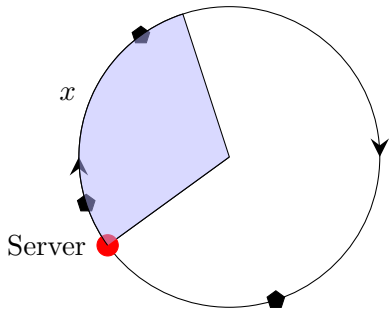

We define the distance between a point on the circle and the server as the distance that the server has to travel to reach that point, see Figure \ref{fig:model}. Throughout this paper we rely on the symmetry of the model, and often instead use that {\em arriving customers obtain locations at a uniform distance from the server}. Additionally we focus on customers at a certain distance to the server, rather than at a location on the circle. We define $L(x)$ to be the number of waiting customers within distance $x$ of the server, see Figure \ref{fig:model}; we often write the shorthand $L$ for all waiting customers in the system (i.e. $L(1)$).

\section{The mean number of waiting customers}
\label{sec:MeanCust}
In this section we derive an expression for the expectation of the steady-state number of waiting customers in the system.
We use an approach based on Mean Value Analysis. In two remarks, we briefly outline two other approaches.

Relating this model and one where customers are served in a FCFS order, we prove the following:
\begin{lemma}
    \label{lemma:TotCust}
    The expected number of waiting customers in the system, excluding the customer that is possibly in service, is given by:
    \begin{align}
    \label{eq:totcust}
        \mathbb{E}[L] = \frac{\lambda\mathbb{E}[K]}{2(1-\rho)}\Big(\alpha + \lambda\mathbb{E}[K]\mathbb{E}[B^2]+\mathbb{E}[B]\frac{\mathbb{E}[K(K-1)]}{\mathbb{E}[K]}\Big).
    \end{align}
\end{lemma}
\begin{proof}
 The number waiting in the system, $L$, can be expressed as $L_0 + L_1$, with $L_0 = L$ if 
the server is traveling and $0$ otherwise, and $L_1 = L$ if the server is busy and $0$ otherwise. According to Little's law, $\mathbb{E}[L_0] = \lambda  \mathbb{E}[K]  \mathbb{E}[W_0]$ and $ \mathbb{E}[L_1] = \lambda \mathbb{E}[K] \mathbb{E}[W_1]$, where $W_0$ is the travel time of the server during the waiting time and $W_1$ is the busy time of the server during the waiting time. Clearly, $\mathbb{E}[W_0] = \alpha/2$ and thus
$\mathbb{E}[L_0] = \lambda \mathbb{E}[K]\alpha / 2$. 
To determine $\mathbb{E}[W_1]$, we interchange locations of customers to preserve service in order of arrival. So when a new batch arrives, these customers swap their locations with waiting customers that are overtaken. Interchanging locations of customers does not alter the number waiting in the system, and hence, by Little's law, it does not alter the mean waiting time, since all customers have the same service time distribution. In case of service in order of arrival, we have
\[
E\mathbb{E}[W_1] = \rho \mathbb{E}[R] + (\mathbb{E}[L_0]+\mathbb{E}[L_1]) \mathbb{E}[B] + \frac{\mathbb{E}[K(K-1)]}{2\mathbb{E}[K]} \mathbb{E}[B].
\]
Combined with $\mathbb{E}[L_1] = \lambda \mathbb{E}[K] \mathbb{E}[W_1]$, both $\mathbb{E}[W_1]$ and $\mathbb{E}[L_1]$ follow.
\end{proof}
\begin{remark}
    The expected number of waiting customers can also be analyzed by first considering a discrete cyclic system with $N$ queues, and then letting $N$ tend to $\infty$;
\cite{Eliazar2005} has proven that such a limiting operation, without batch arrivals, results in a continuous polling system. One can extend this result to the case of batch arrivals.
Using the pseudo-conservation law for a discrete polling system with batch arrivals, cf.\ Equation (3.21) of \cite{Boxma1989}, in combination with the symmetry of the model, then results once more in (\ref{eq:totcust}).
\end{remark}
\begin{remark}
The lemma can also be proven by using the stochastic decomposition in \citet[Theorem 2.1]{Boxma1989}: the expected queue length is given by the sum of the expected queue length in a corresponding system without travel time and the expected queue length at an epoch at which the server is traveling.
The first term follows from a known result for the $M/G/1$ queue with batch arrivals (see, e.g., \cite{Cooper1972}, Section 5.10) and equals the last two components of \eqref{eq:totcust}.
The second term can be argued to equal the mean number of customers {\em at the end of a cycle} (one complete tour along the circle) of the server. A balancing argument shows that the mean number of customers arriving per cycle equals $\lambda \mathbb{E}[K] \alpha/(1-\rho)$, where $\alpha/(1-\rho)$ denotes the mean cycle time.
A symmetry argument finally shows that, on average, half of those customers are still present at the end of the cycle.
\end{remark}

\section{The mean sojourn time of a batch}
\label{sec:BatchSojourn}
In this section we derive the expected sojourn time $\mathbb{E}[S^B]$ of a job/batch, i.e., the time from the arrival of a batch until the service completion
of its last customer. The service of jobs, now, is preempted by travel periods, as the server will have to travel between the customers. Moreover, it might happen that other jobs are (partially) served between the service of two customers in the same batch. Because of that, one cannot apply the stochastic decomposition and pseudo-conservation law  in \citet{Boxma1989} to the number of batches in the system.
Making use of a branching process construction, we analyze the total service time of customers served before the tagged customer, allowing us to derive an initial expression for the mean sojourn time of a batch.

Our approach consists of five steps, which are worked out in detail in the following five subsections.
{\bf Step 1} is the branching process construction.
It results in an initial expression for $\mathbb{E} [S^B]$, that still involves one unknown function $f(\cdot)$.
{\bf Step 2} is the construction of an integral equation for that function $f(\cdot)$.
{\bf Step 3} is the solution of that integral equation; its uniqueness is discussed in {\bf Step 4}.
Finally, $\mathbb{E}[S^B]$ is determined in {\bf Step 5}.

\subsection{Step 1: A branching process construction}
The waiting time of a tagged customer consists of the travel time of the server and the service time of all customers served before her. This latter term includes customers who might arrive later than the tagged customer, but with a location closer to the server. Inspired by \citet{Resing1993}, we handle these future arrivals by linking each future arrival to either (i) a customer already present at the arrival moment of the tagged customer or (ii) the travel period, resulting in a multi-type branching process with immigration. 

 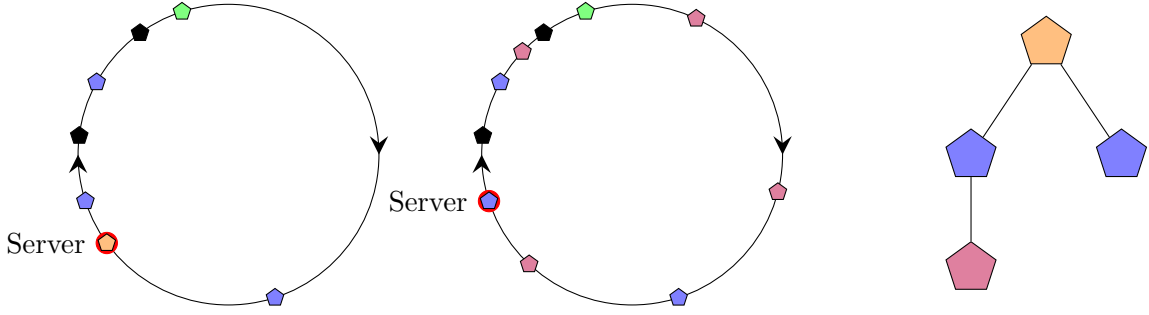
\begin{figure}[!htbp]
    \centering
    \begin{subfigure}{0.3\textwidth}
        \begin{tikzpicture}
        \draw[postaction = {decorate, decoration = {markings, mark = at position 0.6 with {\node[draw, red, circle, fill, inner sep = 1mm] (server){};}}}]
        [postaction = {decorate, decoration = {markings, mark = at position 0.48 with {\node[regular polygon, regular polygon sides = 5,draw, black,  fill, inner sep = 0.7mm] {};}}}]
        [postaction = {decorate, decoration = {markings, mark = at position 0.35 with {\node[regular polygon, regular polygon sides = 5,draw, black,  fill, inner sep = 0.7mm]{};}}}]
        [postaction = {decorate, decoration = {markings, mark = at position 0.6 with {\node[regular polygon, regular polygon sides = 5,draw,  fill=orange!50, inner sep = 0.7mm]{};}}}]
        [postaction = {decorate, decoration = {markings,  mark = at position 0.675 with {\node[inner sep = 0mm] (xlabel) {};}}}]
        [postaction = {decorate, decoration = {markings, mark = at position 0.42 with {\node[regular polygon, regular polygon sides = 5,draw, fill = blue!50,  fill, inner sep = 0.7mm] {};}}}]
        [postaction = {decorate, decoration = {markings, mark = at position 0.55 with {\node[regular polygon, regular polygon sides = 5,draw, fill = blue!50,  fill, inner sep = 0.7mm] {};}}}]
        [postaction = {decorate, decoration = {markings, mark = at position 0.8 with {\node[regular polygon, regular polygon sides = 5,draw, fill = blue!50,  fill, inner sep = 0.7mm] {};}}}]
        [postaction = {decorate, decoration = {markings, mark = at position 0.3 with {\node[regular polygon, regular polygon sides = 5,draw, fill=green!50, inner sep = 0.7mm] {};}}}]
        [postaction = {decorate, decoration = {markings,  mark = at position 0.5 with {\arrowreversed[line width = 0.7mm]{stealth}}}}]
        [postaction = {decorate, decoration = {markings,  mark = at position 1 with {\arrowreversed[line width = 0.7mm]{stealth}}}}]
        (0,0) circle (2);
        \node [left] at (server.west) {Server};
        \end{tikzpicture}
    \end{subfigure}
        \begin{subfigure}{0.3\textwidth}
            \begin{tikzpicture}
            \draw[postaction = {decorate, decoration = {markings, mark = at position 0.55 with {\node[draw, red, circle, fill, inner sep = 1mm] (server){};}}}]
            [postaction = {decorate, decoration = {markings, mark = at position 0.48 with {\node[regular polygon, regular polygon sides = 5,draw, black,  fill, inner sep = 0.7mm] {};}}}]
            [postaction = {decorate, decoration = {markings, mark = at position 0.35 with {\node[regular polygon, regular polygon sides = 5,draw, black,  fill, inner sep = 0.7mm]{};}}}]
            [postaction = {decorate, decoration = {markings,  mark = at position 0.675 with {\node[inner sep = 0mm] (xlabel) {};}}}]
            [postaction = {decorate, decoration = {markings, mark = at position 0.42 with {\node[regular polygon, regular polygon sides = 5,draw, fill = blue!50,  fill, inner sep = 0.7mm] {};}}}]
            [postaction = {decorate, decoration = {markings, mark = at position 0.55 with {\node[regular polygon, regular polygon sides = 5,draw, fill = blue!50,  fill, inner sep = 0.7mm] {};}}}]
            [postaction = {decorate, decoration = {markings, mark = at position 0.8 with {\node[regular polygon, regular polygon sides = 5,draw, fill = blue!50,  fill, inner sep = 0.7mm] {};}}}]
            [postaction = {decorate, decoration = {markings, mark = at position 0.3 with {\node[regular polygon, regular polygon sides = 5,draw, fill=green!50, inner sep = 0.7mm] {};}}}]
            [postaction = {decorate, decoration = {markings,  mark = at position 0.5 with {\arrowreversed[line width = 0.7mm]{stealth}}}}]
            [postaction = {decorate, decoration = {markings,  mark = at position 1 with {\arrowreversed[line width = 0.7mm]{stealth}}}}]
            [postaction = {decorate, decoration = {markings, mark = at position 0.18 with {\node[regular polygon, regular polygon sides = 5,draw, fill=purple!50, inner sep = 0.7mm] {};}}}]
            [postaction = {decorate, decoration = {markings, mark = at position 0.96 with {\node[regular polygon, regular polygon sides = 5,draw, fill=purple!50, inner sep = 0.7mm] {};}}}]
            [postaction = {decorate, decoration = {markings, mark = at position 0.38 with {\node[regular polygon, regular polygon sides = 5,draw, fill=purple!50, inner sep = 0.7mm] {};}}}]
            [postaction = {decorate, decoration = {markings, mark = at position 0.63 with {\node[regular polygon, regular polygon sides = 5,draw, fill=purple!50, inner sep = 0.7mm] {};}}}]
            (0,0) circle (2);
            \node [left] at (server.west) {Server};
            \end{tikzpicture}
    \end{subfigure}
    \hspace{0.04\textwidth}
        \begin{subfigure}{0.3\textwidth}
            \begin{tikzpicture}
            \useasboundingbox (-3,-3.5) rectangle (3,0);
            \node[regular polygon, regular polygon sides = 5, black, draw, fill = orange!50, inner sep = 2mm] (1) at (0,0){};
            \node[regular polygon, regular polygon sides = 5, black, draw, fill = blue!50, inner sep = 2mm] (2) at (-1,-1.5){};
            \node[regular polygon, regular polygon sides = 5, black, draw, fill = blue!50, inner sep = 2mm] (3) at (1,-1.5){};
            \node[regular polygon, regular polygon sides = 5, black, draw, fill = purple!50, inner sep = 2mm] (4) at (-1,-3){};
            \path[draw] (1) -- (2);
            \path[draw] (1) -- (3);
            \path[draw] (2) -- (4);
            \end{tikzpicture}
    \end{subfigure}
    \caption{Illustration of the waiting time of a tagged customer (green) that is generated by a service (of the orange customer) and the corresponding branching process. During the service of the orange customer, blue customers arrive, of which only the first two are considered. During the service of the first blue customer, the red customers arrive, of which only one will be served before the tagged customer.}
    \label{fig:generatedwaiting}
\end{figure}

We construct the branching process by following two simple rules:
\begin{itemize}
    \item A customer (a) is the offspring of a customer (b) when (a) arrives during the service of (b)
    \item A customer (a) is called an immigrant when she arrives during a travel period of the server.
\end{itemize}
From the perspective of a tagged customer, however, only customers who will be served before her are of interest. We therefore trim the branching process to only consider the customers who are served before the tagged customer, see Figure \ref{fig:generatedwaiting}. As the offspring distribution of a customer now depends on the distance from this customer to the tagged customer, we take this distance as the type of a customer. Using this trimmed branching process we introduce the following.
\begin{definition}
$S(y)$ denotes the service time of a customer at distance $y$ of the tagged customer plus the service times of all its descendants in the trimmed branching process. This is also referred to as the
\textit{waiting time generated by a service at distance $y$}.
\end{definition}
\begin{definition}
$T(y)$ denotes the added service times of all customers in the branching processes of immigrants, arriving during the travel time of the server of distance $y$ to the tagged customer. We call this the \textit{waiting time generated by the traveling of distance $y$}.
\end{definition}

The expectations of these random variables can be found by conditioning w.r.t.\ the size of the offspring, resulting in an integral equation for the first moment of $S(y)$,
that allows us to obtain $\mathbb{E}[S(x)]$ and $\mathbb{E}[T(x)]$: 

\begin{lemma}
\label{lemma:generatedwait}
    The expected waiting time generated by a service at distance $x$ of a customer, $\mathbb{E}[S(x)]$, is given by:
    \begin{align}
        \label{eq:servicework_nobatch}
        \mathbb{E}[S(x)] = \mathbb{E}[B]\exp(\rho x).
    \end{align}
    The mean total waiting time generated by the server traveling a distance $x$ to the customer, $\mathbb{E}[T(x)]$, satisfies:
    \begin{align}
        \label{eq:travelwork_nobatch}
        \mathbb{E}[T(x)] = \frac{\alpha}{\rho}\Big[\exp(\rho x) - 1\Big].
    \end{align}
\end{lemma}
\begin{proof}
    During a service, we know that in expectation a total of $\rho$ customers arrive. Each of them initiates another additional service requirement of $S(y)$, where $y$ is the location of the arriving customer. Hence: 
\begin{align*}
    \mathbb{E}[S(x)] = \mathbb{E}[B] + \rho \int_{y= 0}^x \mathbb{E}[S(y)]\dx y \Longrightarrow \mathbb{E}[S(x)] = \mathbb{E}[B]\exp(\rho x),
\end{align*}
where the last step follows from the fact that $\mathbb{E}[S(0)] = \mathbb{E}[B]$.\\
We can now use the above expression to prove (\ref{eq:travelwork_nobatch}). The server takes a time $\alpha \dx u$ to travel a small distance $\dx u$. Therefore, in expectation, $\lambda \mathbb{E}[K]\alpha \dx u$ customers arrive during the traveling of said distance. Each of these arrivals adds a total waiting time of $S(y)$ when  arriving between the tagged customer and the server at a distance $y$ to the tagged customer:
\begin{align*}
    \mathbb{E}[T(x)] &= \alpha x +\int_{u=0}^x \lambda\mathbb{E}[K]\alpha \int_{y=0}^{u} \mathbb{E}[S(y)]\dx y \dx u = \alpha x +\alpha \lambda\mathbb{E}[K] \cdot \int_{u=0}^x \frac{1}{\lambda\mathbb{E}[K]}\Big[\exp(\rho u) - 1\Big]\dx u\\
    &=\frac{\alpha}{\rho}\Big[\exp(\rho x) - 1\Big]. \qedhere
\end{align*}
\end{proof}
\begin{remark}
    At the arrival instant of a customer she might encounter a residual service time first. The total expected waiting time generated by this residual service time, $\mathbb{E}[S^R(x)]$, satisfies relation \eqref{eq:servicework_nobatch} with $\mathbb{E}[B]$ replaced by $\mathbb{E}[B^R]$.
\end{remark}

For a direct analysis of the mean batch sojourn time we first condition on the distance of the furthest customer w.r.t. the server, denoted as $X^B$. This is the maximum of $K$ uniform random variables, hence conditional on $K=k$ the c.d.f. of $X^B$ is given by $x^k$, therefore:
\begin{align}
\label{eq:WB1}
        \mathbb{E}[S^B] = \mathbb{E}[B] + \int_{x=0}^1 \sum_{k=1}^\infty p_kkx^{k-1} \mathbb{E}[W^B\vert X^B = x, K = k]\dx x,
\end{align}
where $W^B$ denotes the waiting time of the last customer of a job. Remark that this waiting time of a customer can be split in four distinct elements: (i) the total waiting time generated by the traveling of the server over distance $x$: $T(x)$; (ii) the waiting time generated by the customer currently in service, if any: with probability $\rho$ there is a contribution $S^R(x)$, where $S^R(x)$ denotes the total waiting time generated by the \emph{residual} service of a customer; (iii) the waiting time generated by all customers in the same batch $S(Y_1) + ... + S(Y_{k-1})$, where $Y_1,...,Y_{k-1}$ denote the distance of these (unordered) customers to the tagged customer -- remark that these are uniform on $[0,x]$; (iv) the waiting time generated by all customers present at the time of arrival and in front of the tagged customer. A customer at distance $y$ 
of the server has a distance $x-y$ to the tagged customer, and hence generates a waiting time of $S(x-y)$. Remembering that $L(x)$ denotes the number of waiting customers within distance $x$ of the server, let $f(x)$ denote the derivative of $\mathbb{E}[L(x)]$. 
Intuitively, $f(x) {\rm d}x$ denotes the long-run expected number of customers with distances between $x$ and $x+\dx x$ of the server.
Then the mean of the maximal waiting time of a batch, conditional on furthest customer location and the batch size, is given by
\begin{align}
\label{eq:WBcond1}
    \mathbb{E}[W^B\vert X^B = x, K = k] &= \mathbb{E}[T(x)] + \rho\mathbb{E}[S^R(x)]
     + (k-1)\int_{y=0}^x \frac{1}{x} \mathbb{E}[S(y)]\dx y\\ &\quad+ \int_{y=0}^x \mathbb{E}[S(x-y)]f(y)\dx y.
     \nonumber
\end{align}

Equation \eqref{eq:WBcond1} reveals that, for a direct analysis of the waiting time of a job/batch,  we require an expression for $f(x)$. Such an expression can be found by constructing an integral equation for $f(x)$, based on Little's law applied to intervals.

\subsection{Step 2: An integral equation}
In this subsection we construct an integral equation that $f(x)$ has to satisfy. The foundation of this subsection is an observation similar to Little's law. We relate the number of customers with distances further than $1-x$ from the server to the waiting time of an arbitrary customer arriving in the range $[0,x]$.

\begin{lemma}[Little's law]
    \label{lemma:RestrictedLittle}
    Let $W(x)$ be the waiting time of an arbitrary customer arriving within distance $x$ of the server. Then we have:
    \begin{align}
    \label{eq:RestrictedLittle}
        \mathbb{E}[L] - \mathbb{E}[L(1-x)] = \lambda\mathbb{E}[K] x\mathbb{E}\big[W(x)\big].
    \end{align}
\end{lemma}
\begin{proof}
    Consider the case in which each customer obtains $\$1$ for each time unit it spends at a distance of at least $1-x$ of the server. One could pay this out by giving each customer the dollar per time unit, i.e. each time unit you pay: $\mathbb{E}[L(1)] - \mathbb{E}[L(1-x)]$ in expectation (observe that $L(1)=L$).\\    
    Equivalently one could pay out each customer leaving the interval $[1-x,1]$, where the customer thus receives a dollar for each time unit she spends in this interval: denoted as $\tau_{[1-x,1]}$. As the departure rate from this interval has to equal the arrival rate, due to the stability condition, this implies that each time unit you expect to pay out $\lambda \mathbb{E}[K] x$ customers.\\
    Since both payout methods are equivalent we thus have:
    \begin{align*}
        \mathbb{E}[L] - \mathbb{E}[L(1-x)] =\lambda \mathbb{E}[K] x \mathbb{E}[\tau_{[1-x,1]}].
    \end{align*}
    Now consider a customer arriving at distance $1-x+y$ and an imaginary customer arriving at $y$. Then the total time the first customer spends in the interval $[1-x,1]$ is equal to the total time it takes the server to reach point $y$, hence the waiting time of the imaginary customer. Therefore, due to the symmetry of the model, we have: $\mathbb{E}[\tau_{[1-x,1]}]= \mathbb{E}\big[W(x)\big]$.
\end{proof}

\begin{remark}
     This observation is similar to that in \citet[Equation (6)] {Winands2006}, regarding discrete polling models. The equivalence becomes especially apparent when one no longer focuses on customers at distance $x$ of the server but rather on the customers at a certain fixed location.
\end{remark}

Lemmas \ref{lemma:generatedwait} and \ref{lemma:RestrictedLittle} form the main ingredients for the integral equation that we construct in this subsection.
We write the expected waiting time $\mathbb{E}[W(x)]$ as an integral by conditioning over the arrival location of the tagged customer, and by then taking the derivative w.r.t. $x$ on both sides of \eqref{eq:RestrictedLittle} we find:
\begin{proposition}
\label{prop:DiffEquation}
The function $f(x)$ satisfies the following integral equation:
\begin{align}
\label{eq:differential1}
   f(1-x) = \lambda\mathbb{E}[K]\bigg\{&\frac{\alpha}{\rho}\big[\exp(\rho x) - 1\big] + \frac{\lambda\mathbb{E}[K] \mathbb{E}[B^2]}{2}\exp(\rho x) + \mathbb{E}[B]\frac{\mathbb{E}[K(K-1)]}{\rho \mathbb{E}[K]}\Big[\exp(\rho x ) - 1\Big]\nonumber\\
        &+\int_{z=0}^x \mathbb{E}[B]\exp\big(\rho(x-z)\big)f(z)\dx z\bigg\}, \quad 0\leq x < 1. 
\end{align}
\end{proposition}
\begin{proof}
First we focus on $\mathbb{E}[W(x)]$ and use that the waiting time consists of the same four parts as before: (i) the waiting time generated by the traveling of the server to the customer; (ii) the waiting time generated by the possible residual service time; (iii) the total waiting time generated by the customers in the same batch who arrive in front of the customer; (iv) the waiting time generated by customers already present on the circle at the moment of arrival. Conditioning on the exact location of the customer shows:
\begin{align}
    \mathbb{E}\big[W(x)\big] &=    \frac{1}{x}\int_{u=0}^{x}\begin{aligned}[t]
        \bigg\{& \underbrace{\mathbb{E}[T(u)]}_{\text{(i)}} +\underbrace{ \rho \mathbb{E}[S^R(u)]}_{\text{(ii)}} + \underbrace{\frac{\mathbb{E}[K(K-1)]}{\mathbb{E}[K]}\int_{z=0}^u \mathbb{E}[S(u-z)]\dx z}_{\text{(iii)}} \\
        &+\underbrace{\int_{z=0}^u \mathbb{E}[S(u-z)]f(z)\dx z}_{\text{(iv)}}
            \bigg\} \dx u, \end{aligned}  \nonumber \\
\intertext{and hence:}
        x\mathbb{E}\big[W(x)\big] &=\int_{u=0}^{x}\begin{aligned}[t]
        \bigg\{& \frac{\alpha}{\rho}\big[\exp(\rho u) - 1\big] + \frac{\lambda\mathbb{E}[K] \mathbb{E}[B^2]}{2}\exp(\rho u) + \mathbb{E}[B]\frac{\mathbb{E}[K(K-1)]}{\rho \mathbb{E}[K]}\Big[\exp(\rho u) - 1\Big]\\
        &+\int_{z=0}^u \mathbb{E}[B]\exp\big(\rho(u-z)\big)f(z)\dx z\bigg\}\dx u.
    \end{aligned}\nonumber
\end{align}
By writing the LHS of \eqref{eq:RestrictedLittle} as an integral over $f(u)$, i.e. $\mathbb{E}[L]-\mathbb{E}[L(1-x)] = \int_{u=1-x}^1 f(u)\dx u $, we obtain the following equality:
\begin{align}
\label{eq:fpe1}
    \int_{u=1-x}^1 f(u)\dx u &=\lambda\mathbb{E}[K]\int_{u=0}^{x}
        \bigg\{ \frac{\alpha}{\rho}\big[\exp(\rho u) - 1\big] + \frac{\lambda\mathbb{E}[K] \mathbb{E}[B^2]}{2}\exp(\rho u) \\
        &\quad + \mathbb{E}[B]\frac{\mathbb{E}[K(K-1)]}{\rho \mathbb{E}[K]}\Big[\exp(\rho u) - 1\Big]
        +\int_{z=0}^u \mathbb{E}[B]\exp\big(\rho(u-z)\big)f(z)\dx z\bigg\}\dx u.\nonumber
\end{align}
Taking the derivative on both sides results in the integral equation:
\begin{align}
\label{eq:diffeq_1}
    f(1-x) = \lambda\mathbb{E}[K]\bigg\{&\frac{\alpha}{\rho}\big[\exp(\rho x) - 1\big] + \frac{\lambda\mathbb{E}[K] \mathbb{E}[B^2]}{2}\exp(\rho x) + \mathbb{E}[B]\frac{\mathbb{E}[K(K-1)]}{\rho \mathbb{E}[K]}\Big[\exp(\rho x ) - 1\Big]\\
        &+\int_{z=0}^x \mathbb{E}[B]\exp\big(\rho(x-z)\big)f(z)\dx z\bigg\}. \qedhere
\end{align}
\end{proof}
\begin{corollary}
    Any twice-differentiable solution $f(x)$ to \eqref{eq:differential1} satisfies $f''(x) = 0$.
\end{corollary}
\begin{proof}
    Taking the derivative on both sides of \eqref{eq:differential1} results in:
    \begin{align*}
        -f'(1-x) = \rho \lambda \mathbb{E}[K] \bigg\{&\frac{\alpha}{\rho}\big[\exp(\rho x) - 1\big] + \frac{\lambda\mathbb{E}[K] \mathbb{E}[B^2]}{2}\exp(\rho x) + \mathbb{E}[B]\frac{\mathbb{E}[K(K-1)]}{\rho \mathbb{E}[K]}\Big[\exp(\rho x ) - 1\Big]\nonumber\\
        &+\int_{z=0}^x \mathbb{E}[B]\exp\big(\rho(x-z)\big)f(z)\dx z\bigg\} + \alpha \lambda \mathbb{E}[K] + \lambda\mathbb{E}[B]\mathbb{E}[K(K-1)] + \rho f(x).
    \end{align*}
    Here, the first part is $\rho f(1-x)$ by \eqref{eq:diffeq_1} and thus:
    \begin{align*}
        -f'(1-x) = \alpha \lambda \mathbb{E}[K] + \lambda\mathbb{E}[B]\mathbb{E}[K(K-1)] + \rho f(1-x) + \rho f(x).
    \end{align*}
    Taking the derivative on both sides another time and using the above relation results in:
    \begin{align*}
        f''(1-x) &= -\rho f'(1-x) + \rho f'(x)= \rho^2 f(1-x) + \rho^2 f(x) - \rho^2 f(x) - \rho^2 f(1-x) = 0. \qedhere
    \end{align*}
\end{proof}
The corollary hints at a linear solution to the integral equation. This can also be argued intuitively.
All customers at a distance $x$ of the server need to have arrived since the last visit to said point. During this time, the server has traversed a total distance of $1-x$. Due to the symmetry of the model, it makes sense that the time the server takes, on average, to traverse this distance is proportional to that distance. 

\subsection{Step 3: Solution of the integral equation}
In this subsection we give a (linear) solution to the previously constructed integral equation with additional condition that the integral over $f(x)$ has to equal $\mathbb{E}[L]$, found in Lemma \ref{lemma:TotCust}. We then prove that this is the unique solution and therefore that this solution gives the expected number of customers in a small interval around $x$. We then derive the expected batch sojourn time, see Theorem \ref{thm:customers_batch}.

\begin{proposition}
\label{prop:diffSolution}
The function:\\
\begin{align}
\label{eq:solf}
    f(x) &= \rho\frac{\lambda\mathbb{E}[K]\mathbb{E}[B^2]}{2\mathbb{E}[B]} + \frac{\alpha \lambda\mathbb{E}[K]}{1-\rho}(1-x) + \rho\frac{\lambda^2\mathbb{E}[K]^2\mathbb{E}[B^2]}{1-\rho}(1-x) + \frac{\rho\mathbb{E}[K(K-1)]}{(1-\rho)\mathbb{E}[K]}(1-x),
\end{align}
is a solution to \eqref{eq:differential1} with boundary condition:
\begin{align}
\label{eq:BoundaryCondition}
    &\int_{x=0}^1 f(x)\dx x =  \frac{\lambda\mathbb{E}[K]}{2(1-\rho)}\Big(\alpha + \lambda\mathbb{E}[K]\mathbb{E}[B^2]+\mathbb{E}[B]\frac{\mathbb{E}[K(K-1)]}{\mathbb{E}[K]}\Big).
\end{align}
\end{proposition}
\begin{proof}
It is readily verified that $f(\cdot)$ satisfies the boundary condition \eqref{eq:BoundaryCondition}. It is thus left to prove that the expression in \eqref{eq:solf} also satisfies equation \eqref{eq:differential1}.\\
For the linearity of the solution to the integral equation, we substitute $f(1-x) = ax+ b$ in \eqref{eq:differential1} and show that \eqref{eq:solf} gives a solution. By partial integration we have that:
\begin{align*}
    \int_{z=0}^x \rho\exp\big(\rho(x-z)\big)f(z)\dx z &= \bigg[-\big\{a(1-z)+b\big\}\cdot\exp\big(\rho(x-z)\big) \bigg]_{z=0}^x -a  \int_{z=0}^x \exp\big(\rho(x-z)\big)\dx z \\
    &=-a(1-x) - b + (a+b)\exp(\rho x) + \frac{a}{\rho}\big[1-\exp(\rho x)\big] \\
    &= ax + \Big(b - \frac{a(1-\rho)}{\rho}\Big)\big[\exp(\rho x)-1\big].
\end{align*}
Remark that $ax$ cancels against the left-hand side in \eqref{eq:differential1}. Further, by taking $b = \lambda^2\mathbb{E}[K]^2\mathbb{E}[B^2]/2$, the only remaining terms in the right-hand side of \eqref{eq:differential1} are $\exp(\rho x) -1$ factors. These disappear when taking $a$ as in \eqref{eq:solf}, i.e.: 
\begin{equation*}
    a = \frac{1}{1-\rho}\bigg(\alpha\lambda\mathbb{E}[K] + \rho\lambda^2\mathbb{E}[K]^2\mathbb{E}[B^2] + \frac{\rho\mathbb{E}[K(K-1)]}{\mathbb{E}[K]}\bigg).\qedhere
\end{equation*}
\end{proof}

\subsection{Step 4: uniqueness of the solution of the integral equation}
In order to prove that the solution above indeed represents the expected number of customers in a small interval around $x$ we prove that the integral equation with boundary conditions has a unique integrable solution.

\begin{proposition}
The integral equation in \eqref{eq:differential1} with boundary condition \eqref{eq:BoundaryCondition} has a unique integrable solution.
\end{proposition}
\begin{proof}
    We first use that the integral equation implies \eqref{eq:fpe1} by integrating from $1-x$ to $1$. Interchanging the order of integration now gives:
\begin{align*}
        \int_{u=1-x}^1 f(u)\dx u &=\lambda\mathbb{E}[K]\int_{u=0}^{x}\bigg\{
        \begin{aligned}[t]
        &\frac{\alpha}{\rho}\big[\exp(\rho u) - 1\big] + \frac{\lambda\mathbb{E}[K] \mathbb{E}[B^2]}{2}\exp(\rho u) \\
        &+ \mathbb{E}[B]\frac{\mathbb{E}[K(K-1)]}{\rho \mathbb{E}[K]}\Big[\exp(\rho u) - 1\Big]\bigg\}\dx u
        \end{aligned}\nonumber\\
        &\quad+\int_{z=0}^x f(z)\int_{u=z}^x \rho\exp\big(\rho(u-z)\big)\dx u\dx z\\
        &=\lambda\mathbb{E}[K]\int_{u=0}^{x}\bigg\{
        \begin{aligned}[t]
        &\frac{\alpha}{\rho}\big[\exp(\rho u) - 1\big] + \frac{\lambda\mathbb{E}[K] \mathbb{E}[B^2]}{2}\exp(\rho u) \\
        &+ \mathbb{E}[B]\frac{\mathbb{E}[K(K-1)]}{\rho \mathbb{E}[K]}\Big[\exp(\rho u) - 1\Big]\bigg\}\dx u
        \end{aligned}\nonumber\\
        &\quad+\int_{z=0}^x f(z)\big[\exp\big(\rho(x-z)\big)-1\big]\dx z.
\end{align*}
Single out the integral $\int_{z=0}^x f(z)\dx z$ and take it to the left-hand side to conclude:
\begin{align}
    \label{eq:fpe1_v2}
    \int_{u=1-x}^1 f(u)\dx u + \int_{u=0}^x f(u)\dx u &= \lambda\mathbb{E}[K]\int_{u=0}^{x}\bigg\{\begin{aligned}[t] &\frac{\alpha}{\rho}\big[\exp(\rho u) - 1\big] + \frac{\lambda\mathbb{E}[K] \mathbb{E}[B^2]}{2}\exp(\rho u) \\
    &+ \mathbb{E}[B]\frac{\mathbb{E}[K(K-1)]}{\rho \mathbb{E}[K]}\Big[\exp(\rho u) - 1\Big]\bigg\}\dx u
    \end{aligned}\nonumber\\
    &\quad + \int_{u=0}^x f(u)\exp\big(\rho(x-u)\big)\dx u.
\end{align}

Now consider two functions $f,g$ satisfying \eqref{eq:differential1}, and let $||\cdot||_{\infty}$ denote the infinity-norm (supremum norm). Since the first terms of \eqref{eq:differential1} are not affected by the choice of $f,g$, we have:
\begin{align*}
    ||f-g||_{\infty} &=  \Big|\Big|\int_{u=0}^x \rho\exp\big(\rho(x-u)\big)(f(u)-g(u))\dx u\Big|\Big|_{\infty}.
\intertext{It follows from \eqref{eq:fpe1_v2} that we can rewrite this as follows:}
    ||f-g||_{\infty}& =\Big|\Big|\rho\int_{u=1-x}^1 (f(u)-g(u))\dx u + \rho\int_{u=0}^x (f(u)-g(u))\dx u\Big|\Big|_{\infty}.
\end{align*}
If $x \geq 1/2$, then we can write the integral as the integral from $0$ to $1$ plus the part that is accounted for twice:
\begin{align*}
    ||f-g||_{\infty} &= \Big|\Big|\rho\int_{u=0}^1 (f(u)-g(u))\dx u + \rho\int_{u=1-x}^x (f(u)-g(u))\dx u\Big|\Big|_{\infty}\\
    &\leq  \Big|\Big|\rho\cdot 0\Big|\Big|_{\infty} +  \rho\int_{u=1-x}^x||f-g||_{\infty}\dx u\\
    & \leq \rho (2x-1)\cdot ||f-g||_{\infty} \leq \rho ||f-g||_{\infty},
\end{align*}
where the first integral disappeared due to the boundary condition in \eqref{eq:BoundaryCondition}. When $x<1/2$ we have:
\begin{align*}
     ||f-g||_{\infty} &\leq \rho\Big|\Big|\int_{u=1-x}^1 (f(u)-g(u))\dx u\Big|\Big|_{\infty}+ \rho\Big|\Big|\int_{u=0}^x (f(u)-g(u))\dx u\Big|\Big|_{\infty}\\
     &\leq \rho\cdot 2x \cdot ||f-g||_{\infty} \leq \rho ||f-g||_{\infty}.
\end{align*}
Since $\rho < 1$ this implies that $||f-g||_{\infty} = 0$, and hence $f\equiv g$.
\end{proof}

\subsection{Step 5: Determination of $\mathbb{E}[S^B]$}
For the maximal waiting time of a batch we use a construction of the waiting time similar to the proof of Proposition \ref{prop:DiffEquation}. Conditioning on the size of the batch and location of the furthest customer we find:

\begin{theorem}
\label{thm:customers_batch}
    The expected batch sojourn time satisfies:
    \begin{equation}
    \label{eq:batchsojourntime}
        \begin{aligned}[c]
        \mathbb{E}[S^B] &= \mathbb{E}[B] + \frac{1}{1-\rho}\cdot\bigg[\alpha + \rho\lambda\mathbb{E}[K]\mathbb{E}[B^2] + \frac{\mathbb{E}[B]\mathbb{E}[K(K-1)]}{\mathbb{E}[K]}\bigg]\mathbb{E}\Big[\frac{K}{K+1}\Big]\\
          &\quad + \frac{\lambda\mathbb{E}[K]\mathbb{E}[B^2]}{2} +  \frac{\mathbb{E}[B]\mathbb{E}[K(K-1)]}{\rho\mathbb{E}[K]} + \frac{1}{\lambda}\Big(\exp(\rho)-1\Big) \\
          &\quad - \Big(\mathbb{E}[B]+\frac{\mathbb{E}[B]\mathbb{E}[K(K-1)]}{\rho\mathbb{E}[K]}\Big)\exp(\rho)   \\ 
          &\quad +\Big(\mathbb{E}[B]\rho + \frac{\mathbb{E}[B]\mathbb{E}[K(K-1)]}{\mathbb{E}[K]}\Big)\int_{x=0}^1 \exp(\rho x) \tilde{K}(x)\dx x.
        \end{aligned}
    \end{equation}
\end{theorem}
\begin{proof}
    We use, again, that the batch sojourn time is given by the waiting time of the last customer in the batch, plus the service time of said customer. This latter term immediately translates to the first term in \eqref{eq:batchsojourntime}. For the waiting time of the last customer we recall \eqref{eq:WBcond1}:
    \begin{align*}
        \mathbb{E}[W^B\vert X^B = x, K = k] &= \mathbb{E}[T(x)] + \rho\mathbb{E}[S^R(x)] + (k-1)\int_{y=0}^x \frac{1}{x} \mathbb{E}[S(y)]\dx y \\
        &\quad + \int_{y=0}^x \mathbb{E}[S(x-y)]f(y)\dx y\\
      &= \underbrace{\frac{\alpha}{\rho}\big[\exp(\rho x) - 1\big]}_{\text{(i)}} + \underbrace{\frac{\lambda\mathbb{E}[K]\mathbb{E}[B^2]}{2}\exp(\rho x)}_{\text{(ii)}}+\underbrace{\frac{(k-1)\mathbb{E}[B]}{\rho x} \big[\exp(\rho x)-1\big]}_{\text{(iii)}}\\ &\quad+ \underbrace{\int_{y=0}^x \mathbb{E}[B]\exp\big(\rho(x-y)\big)f(y)\dx y}_{\text{(iv)}} .
    \end{align*}
    Applying the same partial integration as in the proof of Proposition \ref{prop:diffSolution}, we find:
    \begin{equation*}
        \text{(iv)} = 
        \begin{aligned}[t]
            &- \frac{\alpha}{\rho}\big(\exp(\rho x ) - 1\big) -\frac{\lambda\mathbb{E}[K]\mathbb{E}[B^2]}{2}\exp(\rho x ) - \frac{\mathbb{E}[B]\mathbb{E}[K(K-1)]}{\rho\mathbb{E}[K]}\big(\exp(\rho x ) - 1\big) \\
            & + \frac{\alpha x}{1-\rho}+\frac{\lambda\mathbb{E}[K]\mathbb{E}[B^2]}{2} + \frac{\rho\lambda\mathbb{E}[K]\mathbb{E}[B^2]x}{1-\rho} + \frac{\mathbb{E}[B]\mathbb{E}[K(K-1)]x}{\mathbb{E}[K](1-\rho)}.
        \end{aligned}
    \end{equation*}
    Now (i) and (ii) cancel against the first elements and thus we are left with:
    \begin{align*}
        \mathbb{E}[W^B\vert X^B = x, K = k] &= \frac{\alpha x}{1-\rho}+\frac{\lambda\mathbb{E}[K]\mathbb{E}[B^2]}{2} + \frac{\rho\lambda\mathbb{E}[K]\mathbb{E}[B^2]x}{1-\rho} + \frac{\mathbb{E}[B]\mathbb{E}[K(K-1)]x}{\mathbb{E}[K](1-\rho)} \\
        &\quad + \frac{(k-1)\mathbb{E}[B]}{\rho x} \big[\exp(\rho x)-1\big] - \frac{\mathbb{E}[B]\mathbb{E}[K(K-1)]}{\rho\mathbb{E}[K]}\big[\exp(\rho x ) - 1\big].
    \end{align*}
    We return to Equation \eqref{eq:WB1} and decondition the expression above with respect to the batch size.    Rewriting everything in terms of the probability generating function of $K$ gives:
    \begin{align*}
        \mathbb{E}[W^B] &= \int_{x=0}^1 \sum_{k=1}^\infty kp_kx^{k-1}\cdot 
        \begin{aligned}[t]
        \bigg[&\frac{\alpha x}{1-\rho}+\frac{\lambda\mathbb{E}[K]\mathbb{E}[B^2]}{2} + \frac{\rho\lambda\mathbb{E}[K]\mathbb{E}[B^2]x}{1-\rho} + \frac{\mathbb{E}[B]\mathbb{E}[K(K-1)]x}{\mathbb{E}[K](1-\rho)}\\
        &+ \frac{(k-1)\mathbb{E}[B]}{\rho x} \big[\exp(\rho x)-1\big] - \frac{\mathbb{E}[B]\mathbb{E}[K(K-1)]}{\rho\mathbb{E}[K]}\big[\exp(\rho x ) - 1\big]\bigg]\dx x
        \end{aligned}\\
        &=\int_{x=0}^1 \bigg[\frac{\alpha}{1-\rho} + \frac{\rho\lambda\mathbb{E}[K]\mathbb{E}[B^2]}{1-\rho} + \frac{\mathbb{E}[B]\mathbb{E}[K(K-1)]}{\mathbb{E}[K](1-\rho)}\bigg]\tilde{K}'(x)x \dx x\\
        &\quad + \int_{x=0}^1 \frac{\lambda\mathbb{E}[K]\mathbb{E}[B^2]}{2}\tilde{K}'(x)\dx x +  \int_{x=0}^1 \frac{\mathbb{E}[B]}{\rho}\big[\exp(\rho x) - 1\big]\tilde{K}''(x) \dx x\\
        &\quad -  \int_{x=0}^1 \frac{\mathbb{E}[B]\mathbb{E}[K(K-1)]}{\rho\mathbb{E}[K]}\big[\exp(\rho x ) - 1\big] \tilde{K}'(x)\dx x.
    \end{align*}
    We now apply partial integration to the first, third and last terms:
    \begin{alignat*}{2}
        &\int_{x=0}^1\tilde{K}'(x)x \dx x &&= \Big[x\tilde{K}(x)\Big]_{x=0}^1 - \int_{x=0}^1\tilde{K}(x) \dx x = \mathbb{E}\Big[\frac{K}{K+1}\Big];\\
         &\int_{x=0}^1\tilde{K}''(x)\big[\exp(\rho x ) - 1\big] \dx x &&= \Big[\big[\exp(\rho x ) - 1\big]\tilde{K}'(x)\Big]_{x=0}^1 - \int_{x=0}^1 \rho \exp(\rho x)\tilde{K}'(x)\dx x \\
         & &&= \big[\exp(\rho) - 1\big]\mathbb{E}[K] - \int_{x=0}^1 \rho \exp(\rho x)\tilde{K}'(x)\dx x; \\
         &\int_{x=0}^1\tilde{K}'(x)\exp(\rho x )\dx x&&= \exp(\rho) -\int_{x=0}^1 \rho\exp(\rho x) \tilde{K}(x)\dx x. 
    \end{alignat*}
    Substituting this in the expression for $\mathbb{E}[W^B]$ results in:
    \begin{align*}
          \mathbb{E}[W^B] &= \bigg[\frac{\alpha}{1-\rho} + \frac{\rho\lambda\mathbb{E}[K]\mathbb{E}[B^2]}{1-\rho} + \frac{\mathbb{E}[B]\mathbb{E}[K(K-1)]}{\mathbb{E}[K](1-\rho)}\bigg]\mathbb{E}\Big[\frac{K}{K+1}\Big]\\
          &\quad + \frac{\lambda\mathbb{E}[K]\mathbb{E}[B^2]}{2} + \frac{\mathbb{E}[B]}{\rho}\Big(\exp(\rho)\mathbb{E}[K] - \mathbb{E}[K] - \rho\exp(\rho) + \rho\int_{x=0}^1 \rho\exp(\rho x) \tilde{K}(x)\dx x\Big)\\
          &\quad - \frac{\mathbb{E}[B]\mathbb{E}[K(K-1)]}{\rho\mathbb{E}[K]}\Big(\exp(\rho) - 1 -\int_{x=0}^1 \rho\exp(\rho x) \tilde{K}(x)\dx x\Big). \qedhere
    \end{align*}
\end{proof}

\begin{remark}
    Taking $\mathbb{E}[B] \to 0$, gives a total mean batch sojourn time of $\alpha\cdot \mathbb{E}[K/(K+1)]$. Note that this is exactly equal to the expected travel time of the server to the furthest customer.
\end{remark}
\begin{remark}
    The case $K \equiv 1$, i.e. unit batch sizes, results in:
    \begin{align*}
        \mathbb{E}[S^B] = \mathbb{E}[B] + \frac{\lambda \mathbb{E}[B^2]}{2(1-\rho)} + \frac{\alpha}{2(1-\rho)}.
    \end{align*}
    This is in line with expectations as the system now     simplifies to the system discussed in \citet{Kroese1992}. 
\end{remark}
\begin{remark}
\label{rem:lambdazero}
    The batch sojourn time for light traffic, i.e. $\lambda \to 0$, can be found by reordering the terms:
        \begin{equation*}
    \label{eq:batchsojourntime}
        \begin{aligned}[c]
        \mathbb{E}[S^B] &= \frac{1}{1-\rho}\cdot\bigg[\alpha + \rho\lambda\mathbb{E}[K]\mathbb{E}[B^2] + \frac{\mathbb{E}[B]\mathbb{E}[K(K-1)]}{\mathbb{E}[K]}\bigg]\mathbb{E}\Big[\frac{K}{K+1}\Big]\\
          &\quad + \frac{\lambda\mathbb{E}[K]\mathbb{E}[B^2]}{2} + \frac{1}{\lambda}\Big(\exp(\rho)-1\Big) \\
          &\quad + \Big(\rho\mathbb{E}[B]+\frac{\mathbb{E}[B]\mathbb{E}[K(K-1)]}{\mathbb{E}[K]}\Big)\int_{x=0}^1 \exp(\rho x) \big[\tilde{K}(x)-1\big]\dx x.
        \end{aligned}
    \end{equation*}
    Note that the last integral simplifies to $-\mathbb{E}[K/(K+1)]$ as $\rho \to 0$. Further realizing that $[\exp(\lambda\mathbb{E}[K]\mathbb{E}[B])-1]/\lambda \to \mathbb{E}[K]\mathbb{E}[B]$ as $\lambda \to 0$, then shows that:
    \begin{align*}
        \mathbb{E}[S^B] \to \alpha\mathbb{E}\Big[\frac{K}{K+1}\Big] + \mathbb{E}[K]\mathbb{E}[B],
    \end{align*}
    i.e. in light traffic the mean batch sojourn time reduces to the sum of the mean travel time to the furthest customer plus the mean service time of a batch.
\end{remark}

\section{Numerical results}
\label{sec:NumericalResults}
The results in this paper allow for a direct performance analysis and can be used to obtain direct insights in the performance effect of, for instance, the distribution of the batch size. Additionally, the given expression lends itself nicely for optimization.  In previous work, \citet{Gaast2017} present a mean value analysis for the discrete polling system with batch arrivals.
They derive the expected waiting time of a batch by first solving a system of $N^2$ equations, in which $N$ is the total number of queues. For symmetric systems, this is reduced to a system of $N$ linear equations. For systems with large $N$ this becomes computationally expensive, and in particular insights into the effect of e.g.\ the distribution of the batches, can no longer be directly obtained. In this section, we approximate a discrete (symmetric) polling system by its continuous variant and investigate the rate of convergence of the mean batch sojourn time as well as the effect of the batch size distribution. 

\begin{figure}[H]
    \centering
    \includegraphics[width = \textwidth]{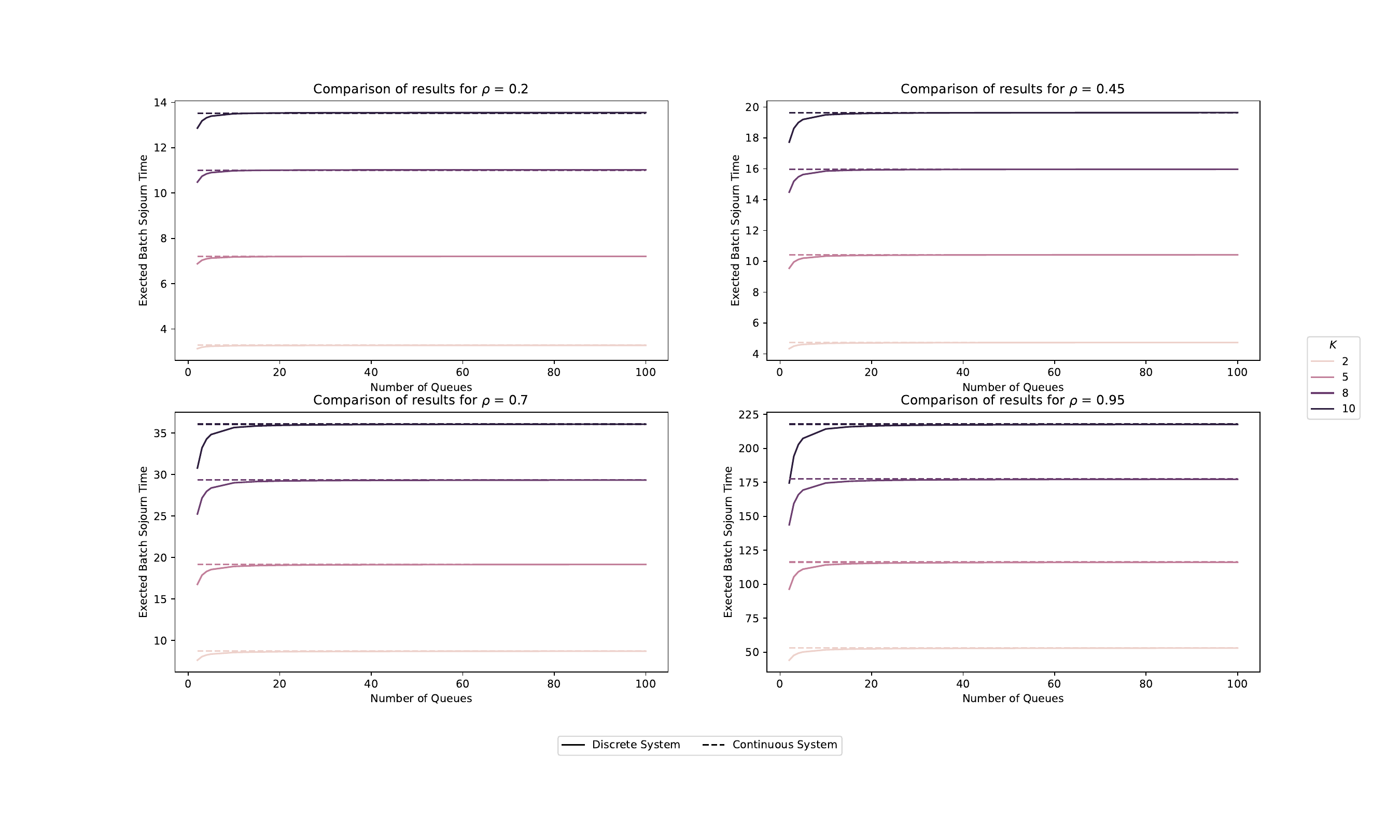}
    \caption{The expected batch sojourn time for deterministic batch sizes and exponential service requirements with unit mean, comparing the discrete and continuous polling model.}
    \label{fig:comp1}
\end{figure}

\begin{figure}[H]
    \centering
    \includegraphics[width = \textwidth]{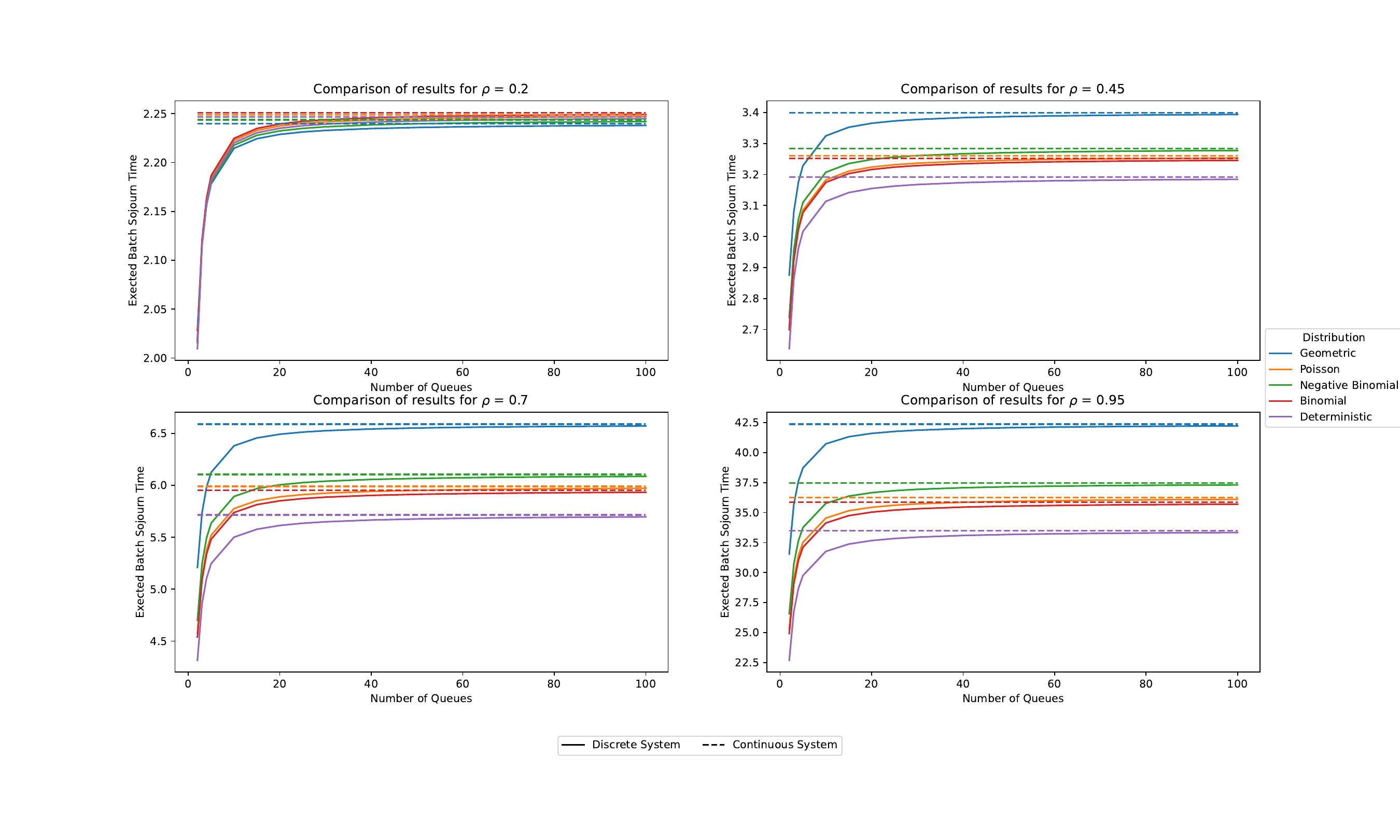}
    \caption{Comparison of the expected batch sojourn time for different batch size distributions: geometric, Poisson, negative binomial with 5 successes, binomial with 15 trials and deterministic. Services take $1/5$ time units and the average batch size is taken $5$.}
    \label{fig:comp2}
\end{figure}

We consider symmetric discrete polling systems, consisting of $N$ queues and deterministic equal switch-over times, with a total switch-over time of $1$. Additionally, we assume that each customer in a batch independently is assigned a queue with equal probability ($1/N$). The generic service time is $B$, the total batch size is $K$. For the continuous variant we take $\alpha$ to be the total switch-over time between queues, i.e. $\alpha = 1$.\\
First, we study the convergence of the mean batch sojourn time for different values of $\rho$ and $K$. Consider the model with deterministic batch sizes and exponential service times (with mean $1$). We consider four values of the workload, $\rho = 0.2,0.45,0.7,0.95$. Figure \ref{fig:comp1} shows the effect of the batch size on the batch sojourn time and illustrates the rate of convergence of the discretized model to the continuous model. Logically, larger batch sizes result in bigger mean sojourn times due to the higher service requirement of a batch. Moreover, the continuous model already proves to be a good approximation for a relatively small number of queues: for $N =  10$ the difference is $<2\%$ and for $N = 20$ the difference is $<1\%$. The rate of the convergence seems to be much higher for lightly loaded systems.\\
Secondly, we consider different batch size distributions, see Figure \ref{fig:comp2}. In each example, the service of a customer is assumed to take a deterministic time with unit mean. Again, it is apparent that, for already small values of $N$, the continuous model provides a reasonable approximation for the batch sojourn time in discrete models, here the difference is $<6\% $ for $N = 10$ and $<3\%$ for $N = 20$. The effect of the batch size distribution on the sojourn time appears to be closely linked to the variance of this distribution; distributions with higher variance result in bigger mean sojourn times. This is in line with results on $M^K/G/1$ queues, where the same relation to the variance of the bulk sizes holds. 

Remark that the continuous system always results in larger mean batch sojourn times. This is caused by the fact that the distance from the server to the furthest customer in a batch is larger in expectation. In a light-traffic situation, i.e. $\lambda \to 0$, this becomes even more apparent as the batch sojourn time is then solely built up from the total service time of a batch and the travel time of the server, cf. Remark \ref{rem:lambdazero}.

\subsubsection*{Acknowledgement}
The research of Tim Engels and Onno Boxma is partly funded by the NWO Gravitation project NETWORKS, grant number 024.002.003.


\bibliography{bib}

\end{document}